\documentclass[11pt]{article}
\usepackage{amssymb}\usepackage{amsmath}
\usepackage{amsthm}
\newcommand{\tr}{\mathop{\rm tr}\nolimits}
\def\ST{{\mathsf T}}

\def\SR{{\mathsf R}}

\def\CT{{\mathcal T}}

\newtheorem{theorem}{\bfseries Theorem}
\newtheorem{propnn}{Proposition}
\newtheorem{lemma}[theorem]{Lemma}
\begin{document} 
 
\begin{center}
{\bf\large\sc   
On orthogonal projections related to representations of the Hecke algebra on a tensor space
} 

\vspace*{3mm}
{ Andrei Bytsko   }
\vspace*{2mm}

\end{center}
\vspace*{2mm}
 
\begin{abstract} 
We consider the problem of finding orthogonal projections $P$ of a rank $r$ that  
give rise to representations of the Hecke algebra $H_N(q)$ in which 
the generators of the algebra act locally on the $N$-th tensor 
power of the space ${\mathbb C}^n$. It is shown that such projections 
are global minima of a certain functional. It is also shown that 
a characteristic property of such projections  is that a certain 
positive definite matrix $A$ has only two eigenvalues or only 
one eigenvalue if $P$ gives rise to a representation of the 
Temperley-Lieb algebra.   
Apart from the parameters $n$, $r$, and $Q=q + q^{-1}$, an additional 
parameter $k$ proves to be a useful characteristic of a projection $P$. 
In particular, we use it to provide a lower bound for $Q$ 
when the values of $n$ and $r$ are fixed and we show that $k= r n$ 
if and only if $P$ is of the Temperley-Lieb type. Besides, we propose 
an approach to constructing projections $P$ and give some novel examples 
for $n=3$.  
\end{abstract} 
 
\section{Introduction}

The Hecke algebra $H_N(q)$, $q \in {\mathbb C}^*$, is the  
unital associative algebra over ${\mathbb C}$ with  
generators $\SR_1, \ldots, \SR_{N-1}$ satisfying the following relations 
\begin{align}
\label{defH1} 
{}& \SR_i^2 = {\mathsf 1} + (q - q^{-1}) \, \SR_i  , \qquad i =1,\ldots, N-1 , \\[1mm]
\label{defH2} 
{}&   \SR_i \SR_m \SR_i = \SR_m \SR_i \SR_m , \qquad  |i-m| =1 ,\\[1mm] 
\label{defH3} 
{}&  \SR_i \SR_m =  \SR_m \SR_i ,  \qquad  |i - m| \geq 2 . 
\end{align}

Let $Q$ be a real positive number.  
Consider a matrix $P \in \mathrm {Mat}(n^2,{\mathbb C})$, $n \geq 2$, that 
satisfies the following relations   
\begin{align} 
{}& \label{PPPm} 
	P^*= P , \qquad  P^2 = P  , \qquad 
 Q^2 ( P_{1} \, P_{2} \, P_{1} - P_{2} \, P_{1} \, P_{2} )
  = P_{1}  - P_{2}  ,     
\end{align}
where $P_{1} \equiv P \,{\otimes}\, I_n$,  
$P_{2} \equiv I_n \,{\otimes}\, P$, $\otimes$ denotes the Kronecker product, 
$I_n \in \mathrm {Mat}(n,{\mathbb C})$ is the identity matrix 
and $P^*$ stands for the  conjugate transpose of $P$. 
The first two relations in (\ref{PPPm}) imply that $P$ is an orthogonal projection 
or, equivalently, a Hermitian idempotent. 
The homomorphism   
$\tau: H_N(q) \to \mathrm {End}\bigl({\mathbb C}^n \bigr)^{\otimes N}$ such that  
\begin{align}\label{tauH} 
 \tau(\SR_i) =  q \,I_n^{\otimes N} - 
 Q \, I_n^{\otimes (i-1)} \otimes P \otimes I_n^{\otimes (N-i-1)}  
\end{align}
is a representation of the algebra $H_N(q)$, where $q+ q^{-1}=Q$. 
In this representation,   $ \tau(\SR_i)$ is Hermitian if $Q \geq 2$ and unitary 
if $Q \leq 2$.
 
Equations (\ref{PPPm}) admit two trivial solutions: 
$P=0$ and $P=I_n \,{\otimes}\, I_n$ (in both 
these cases, $Q$ can be arbitrary). We will call other solutions nontrivial. 
We will say that a solution to (\ref{PPPm}) is of the Temperley-Lieb type if 
$P$ satisfies relations   
\begin{align} 
{}& \label{TTTm} 
	P^*= P , \qquad  P^2 = P  , \qquad 
 Q^2   P_{1} \, P_{2} \, P_{1} =P_{1} , \qquad  
 Q^2  P_{2} \, P_{1} \, P_{2}   = P_{2}  .     
\end{align}
For such a solution $P$, the homomorphism $ \tilde{\tau}(\ST_i) =  
 Q \, I_n^{\otimes (i-1)} \otimes P \otimes I_n^{\otimes (N-i-1)}$ 
provides a representation of the Temperley-Lieb algebra $TL_N(q)$ 
presented by generators 
$\ST_1,\,{\ldots}\,,\ST_{N-1}$ and relations 
\begin{align}
\label{TL1q}
{}& \ST_i^2 = Q \, \ST_i , \qquad i =1,\ldots, N-1 , \\[1mm]
 \label{TL1}
{}& \ST_i \, \ST_{m}  \ST_i = \ST_i   , \qquad  |i -m| =1 ,\\[1mm] 
\label{TL2}
{}& \ST_i \, \ST_m = \ST_m  \ST_i  ,  \qquad  |i-m| \geq 2 . 
\end{align}

A solution $P$ to equations (\ref{PPPm}) yields via the representation (\ref{tauH}) 
a solution to equation~(\ref{defH2}). That is, $R =q I_n^{\otimes 2} - QP$ 
is a solution to the famous Yang-Baxter equation which is known to play 
an important role in statistical mechanics, integrable systems, quantum groups, 
and knot theory. Moreover, thanks to relation (\ref{defH1}), 
$R(\lambda) = \lambda R - \lambda^{-1} R^{-1}$ is a solution to 
the Yang-Baxter equation with a spectral parameter (see e.g.\ section 3.8 in \cite{Is1}). 
One of the motivations to study equations (\ref{PPPm}) 
is that they are a source of R-matrices of this special type.

Given a solution $P \in \mathrm {Mat}(n^2,{\mathbb C})$ to equations (\ref{PPPm}), 
let us set $K_P = P_1 -P_2 \in \mathrm {Mat}(n^3,{\mathbb C})$ and define 
\begin{align}
\label{Prk1}
  r =  \mathrm {rank}\, P , \qquad k =  \frac{1}{2} \, \mathrm {rank }\, K_P .
\end{align} 
Thus, $n$, $r$, $k$, and $Q$ are parameters of a projection $P$. 
The values of these parameters are not entirely independent. 
Some relations between them will be given below. 

\section{Main results}

This Section contains the main results of the present article. 
The proofs of all statements are given in the Appendix A.

\subsection{Solutions as minima of a functional}

Let $\tr_N$ denote the matrix trace on $\mathrm {Mat}(n^N,{\mathbb C})$. 
In particular, we have 
$\tr_3 I_n^{{\otimes}\,3} = n^3$ and $\tr_3 P_1 =  \tr_3 P_2 = r n$. 

\begin{lemma}\label{trPPrn}
Let $P \in \mathrm {Mat}(n^2,{\mathbb C})$ be an orthogonal projection  
of rank $r$. Then  \\[0.5mm] 
a) $\tr_3 (P_1 P_2) \leq r n$ and the equality occurs if and only if 
$P=0$ or $P=I_n \,{\otimes}\, I_n$. \\[0.5mm]
b) $\tr_3 (P_1 P_2)^2 \leq \tr_3 (P_1 P_2)$ 
and the equality occurs if and only if $P_1 P_2 = P_2 P_1$.
\end{lemma}

(To avoid possible confusion, we remark that here and below, 
$\tr (X)^m$ stands for the trace of the $m$-th power of a matrix~$X$.)  

\begin{lemma}\label{Fpos}
Let $F$ be the following  functional on the subset in 
$\mathrm {Mat}(n^2,{\mathbb C})$ consisting of matrices of rank $r$:
\begin{align}\label{trcond0}
F[P]= \bigl( r n - \tr_3 (P_1 P_2) \bigr)  \, 
\bigl( \tr_3 (P_1 P_2)^2 - \tr_3 (P_1 P_2)^3 \bigr) -
\bigl(  \tr_3 (P_1 P_2) - \tr_3 (P_1 P_2)^2  \bigr)^2  .
\end{align}
If $P \in \mathrm {Mat}(n^2,{\mathbb C})$ is an orthogonal projection  
of rank $r$, then $F[P] \geq 0$. 
\end{lemma}

Let us remark that $F[P]=0$ if $P$ is an idempotent matrix such that 
$P_1 P_2 = P_2 P_1$. However, such $P$ satisfies 
relations (\ref{PPPm}) only if $P_1 = P_2$, that is if $P$ is a trivial solution 
to (\ref{PPPm}) (cf.\,the proof of Lemma~\ref{trPPrn}). 

\begin{propnn}\label{tracecond}
Let $P \in \mathrm {Mat}(n^2,{\mathbb C})$  be an orthogonal projection  
of rank $r$. \\[0.5mm]
a) If $P$ satisfies relations (\ref{PPPm}), then 
\begin{align}\label{trcond1}
 F[P] =0  .
\end{align}
b) If condition (\ref{trcond1}) holds and $P_1 P_2 \neq P_2 P_1$, then 
$P$ satisfies relations (\ref{PPPm}), where 
\begin{align}\label{alphaQ}
   Q^2  = 
 \frac{r n - \tr_3 (P_1 P_2) }{\tr_3 (P_1 P_2) - \tr_3 (P_1 P_2)^2 }.
\end{align}
\end{propnn}

Thus, for the algebra $H_N(q)$, the problem of finding representations 
(\ref{tauH}), where $P$ is an orthogonal projection, is equivalent to the problem 
of finding global minima of the functional $F$ on the set of Hermitian idempotents.   
For the Temperley-Lieb algebra $TL_N(q)$, an analogous statement holds  
for the functional 
$\tilde{F}[P]= r n   \tr_3 (P_1 P_2)^2 - \bigl( \tr_3 (P_1 P_2) \bigr)^2$
(see Theorem~1 in \cite{By1}).

\subsection{On the rank of $(P_1 - P_2)$}

For every solution $P$ to equations (\ref{PPPm}), there exists another solution $\tilde{P}$ 
which is dual to $P$ in the following sense.  

\begin{lemma}\label{dualP} 
If $P \in \mathrm {Mat}(n^2,{\mathbb C})$ is a solution 
to (\ref{PPPm}) of rank $r$, then $\tilde{P} = I_n \,{\otimes}\, I_n -P$ 
is a solution to (\ref{PPPm}) of rank $\tilde{r}= (n^2 -r)$ for the same value of~$Q$.
\end{lemma}

Let us stress that if $P$ is a solution to  (\ref{PPPm}) of the 
Temperley-Lieb type, then its dual $\tilde{P}$  is not of the 
Temperley-Lieb type (unless $P$ is a trivial solution or a special 
solution of rank $r = n^2/2$, see Theorem~1 in~\cite{By3} and 
Proposition \ref{nonTL}b below).

Given a matrix $P \in \mathrm {Mat}(n^2,{\mathbb C})$, 
let us consider the matrix $K_P = P_1 -P_2 \in \mathrm {Mat}(n^3,{\mathbb C})$ 
and define $k$ as in (\ref{Prk1}). 
Notice that $K_P=0$ if and only if $P=\text{const}\, I_n \,{\otimes}\, I_n$. 
Notice also that $K_{\tilde{P}} = - K_P$ and so the parameter $k$ 
is the same for $P$ and~$\tilde{P}$.

\begin{lemma}\label{spekK} 
Let $P \in \mathrm {Mat}(n^2,{\mathbb C})$  be an orthogonal projection. 
Then it is a nontrivial solution 
to (\ref{PPPm}) if and only if $Q > 1$ and the only non-zero eigenvalues of $K_P$ are  
$\lambda_\pm =  \pm \sqrt{1 - Q^{-2}}$.  
\end{lemma}

\begin{propnn}\label{rangeK} 
If $P \in \mathrm {Mat}(n^2,{\mathbb C})$ is a nontrivial solution to (\ref{PPPm}) 
of rank $r$, then the eigenvalues 
$\lambda_+$ and $\lambda_-$ of $K_P$ have the same multiplicity~$k$  in the range 
\begin{align}\label{kbounds}
\frac{1}{2n} r (n^2-r)   \leq k \leq \min (r n , n^3 - r n) . 
\end{align}  
\end{propnn}

Let us remark that for some range of $r$, the lower bound for $k$ can be improved (see 
Section~\ref{ontrQ}).

\begin{propnn}\label{TRsol} 
a) If $P \in \mathrm {Mat}(n^2,{\mathbb C})$ is a nontrivial solution to (\ref{PPPm}) 
of rank $r$, then the following equality holds for all integers $m \geq 1$ 
\begin{align}
\label{trPPmk}
{}& \tr_3 (P_{1} P_{2})^m = r n +  (Q^{-2m} -1) \, k  .
\end{align}
b) Let $P \in \mathrm {Mat}(n^2,{\mathbb C})$  be an orthogonal projection  
of rank $r$. If there exist constants $Q > 1$ and $k \in {\mathbb N}$ such that 
relation (\ref{trPPmk}) holds for $m=1,2,3$, then $P$ is a nontrivial solution to~(\ref{PPPm}).
\end{propnn}

Equation (\ref{trPPmk}) shows that $k = \frac{1}{2} \mathrm {rank }\, K_P$ 
is an important parameter of a solution to (\ref{PPPm}) and thus of the representation~(\ref{tauH}).   
In particular, this parameter characterizes the degree of ``Temperley-Lieb-ness'' 
of a given solution in the following sense.

\begin{propnn}\label{nonTL}
Let $P \in \mathrm {Mat}(n^2,{\mathbb C})$ be a nontrivial solution 
to (\ref{PPPm}) of rank $r$. \\[1mm]
a) The following equality holds 
\begin{align}\label{rankPPP}
 \mathrm {rank }\,  \bigl( Q^2   P_{1} P_{2} P_{1} - P_{1} \bigr) 
 =  \mathrm {rank }\,  \bigl( Q^2   P_{2} P_{1} P_{2} - P_{2} \bigr)
  =r n - k .
\end{align}
b) $P$ is of the Temperley-Lieb type, i.e.~it is a solution to (\ref{TTTm}),
 if and only if $k = r n$.
\end{propnn}

Let us remark that $\Pi_3 = (Q^2-1)^{-1} 
 \bigl( Q^2   \tilde{P}_{1} \tilde{P}_{2} \tilde{P}_{1} - \tilde{P}_{1} \bigr)$, 
where $\tilde{P}$ was defined in Lemma~\ref{dualP}, 
is an orthogonal projection such that $P_1 \Pi_3 = P_2  \Pi_3 =0$. 
Relation (\ref{rankPPP}) applied to $\tilde{P}$, yields 
$\tr_3 \Pi_3   = \tilde{r} n - k = n^3 -r n - k$.

\subsection{On the spectrum of $A_\CT$}

Let $\CT = \{ V_1, \ldots, V_r \}$ 
be a set of $r$ complex matrices of size $n$ such that 
\begin{equation}\label{vv}
   \tr \bigl( V_s^*  V_m \bigr)   =\delta_{sm}  .
\end{equation}
Let $E^{(n)}_{ab} \in \mathrm {Mat}(n,{\mathbb C})$ denote the matrix unit such that 
$(E^{(n)}_{ab})_{ij} = \delta_{ai}  \delta_{bj}$. Then 
\begin{equation}\label{PiTau}
  P_\CT = 
 \sum_{s=1}^r \sum_{a,b,c,d=1}^n (V_s)_{ab} \, 
 (\bar{V}_s)_{cd} \  E^{(n)}_{ac} \otimes E^{(n)}_{bd} ,
\end{equation}
where $\bar{V}$ denotes the complex conjugate of $V$, 
is an orthogonal projection of rank~$r$. 
The condition that $P_\CT$ is a solution to (\ref{PPPm}) can be 
reformulated as a condition on the matrices $V_1, \ldots, V_r$.
To this end, we define the  matrices 
$W_{\CT}, A_\CT \in \mathrm {Mat}(r n,{\mathbb C})$ as follows 
\begin{align}\label{WVAV} 
  W_{\CT}  
  = \sum_{s,m=1}^r  E^{(r)}_{sm} \otimes V_m \bar{V}_s , \qquad
  A_\CT =  W_\CT W_\CT^* =
   \sum_{s,m=1}^r  E^{(r)}_{sm} \otimes
   \Bigl( \sum_{i=1}^r V_i \bar{V}_s V^t_m V_i^* \Bigr) . 
\end{align}

\begin{lemma}\label{trAP} 
Let $P_\CT $ be given by (\ref{PiTau}),
where $V_1, \ldots, V_r$ satisfy condition~(\ref{vv}), 
and let $A_\CT $ be given by (\ref{WVAV}). 
Then the following equality holds for all integers $m \geq 1$ 
\begin{align}
\label{trPPmkA}
{}& \tr_3 \bigl( (P_\CT)_{1} (P_\CT)_{2} \bigr)^m = \tr (A_\CT)^m  .
\end{align}
\end{lemma}

An immediate corollary of this Lemma is the following counterpart 
of Proposition~\ref{TRsol}. 

\begin{propnn}\label{TRAsol} 
Let $P_\CT $ be given by (\ref{PiTau}),
where $V_1, \ldots, V_r$ satisfy condition~(\ref{vv}),  
and let $A_\CT $ be given by (\ref{WVAV}). \\[0.5mm]
a) If $P_\CT$ is a nontrivial solution to (\ref{PPPm}), 
then the following equality holds for all integers $m \geq 1$ 
\begin{align}
\label{trAmk}
{}& \tr (A_\CT)^m = r n +  (Q^{-2m} -1) \, k  .
\end{align}
b)   If there exist constants $Q > 1$ and $k \in {\mathbb N}$ such that 
relation (\ref{trAmk}) holds for $m=1,2,3$, then $P_\CT$ is a nontrivial solution to~(\ref{PPPm}).
\end{propnn}

Since $A_\CT$ is a Hermitian matrix, equation (\ref{trAmk}) implies that 
$A_\CT$ has eigenvalue 1 with multiplicity $(r n -k)$ and 
eigenvalue $Q^{-2}$ with multiplicity~$k$. Taking into account that 
$A_\CT$ is a matrix of size $r n$, we conclude that $A_\CT$ is a 
positive definite matrix whose spectrum contains only two points 
(unless $Q=1$ or $k= r n$).  Thus, together with Proposition~\ref{nonTL}b, 
Proposition~\ref{TRAsol}  implies the following statement (its first part recovers  
Theorem~2 in \cite{By1}).  

\begin{propnn}\label{AEsol} 
Let $P_\CT $ be given by (\ref{PiTau}),
where $V_1, \ldots, V_r$ satisfy condition~(\ref{vv}), 
and let $W_\CT$, $A_\CT $ be given by (\ref{WVAV}). \\[0.5mm]
a) $P_\CT$ is a nonzero solution to (\ref{TTTm}) if and only if 
$Q W_\CT$ is a unitary matrix or, equivalently, if and only if 
$Q^2 A_\CT$ is the identity matrix. \\[1mm]
b) $P_\CT$ is a nontrivial solution to (\ref{PPPm}) not of the Temperley-Lieb type
if and only if the only eigenvalues of $A_\CT$ are 1 and $Q^{-2}$, where $Q>1$. 
\end{propnn}

\subsection{On the range of $Q$}\label{ontrQ}

It is an interesting problem to determine the range of possible values of $Q$ 
for a solution to (\ref{PPPm}) for the given values of $n$ and~$r$.  
In the Temperley-Lieb case, $Q \in \{1, \sqrt{2}, \sqrt{3} \}$ if $r> n^2/4$ 
and $Q \geq \max (2,n/r)$  if $r \leq n^2/4$ (see 
Theorem~3 in~\cite{By1} and Theorem~1 in~\cite{By3}). 

\begin{lemma}\label{QNrk1}
If $P \in \mathrm {Mat}(n^2,{\mathbb C})$ is a solution 
to (\ref{PPPm}) of rank $r$, then the following inequality holds
\begin{align}\label{Qrn0}
     r n -k +  Q^{-1} k \leq r^2. 
\end{align}
\end{lemma}

As a corollary to this Lemma, we obtain the following statement 
(note that Proposition~4.1 in~\cite{Gur1} contains an equivalent assertion).

\begin{propnn}\label{TLrank1} 
Every solution to (\ref{PPPm}) of rank $1$ 
is of the Temperley-Lieb type, i.e.~it is a solution to~(\ref{TTTm}).
\end{propnn}

For a solution $P$ of the Temperley-Lieb type, inequality (\ref{Qrn0}) 
yields the estimate $Q \geq n/r$ 
since we have $k=r n$ by Proposition~\ref{nonTL}b. 
In the general case, we have the following statement. 
 
\begin{propnn}\label{Qmin} 
Let $P \in \mathrm {Mat}(n^2,{\mathbb C})$ be a nontrivial 
solution to (\ref{PPPm}) of rank $r$ not of the Temperley-Lieb type. \\[0.5mm]
a) If $1 < r <n$, then 
\begin{align}
\label{kQmin1}
&{}   r n - r^2 +1 \leq  k  \leq r n - 1, \qquad  Q   \geq  \frac{r n -1}{r^2-1} . 
\end{align}
b) If $n^2 - n < r <n^2$, then 
\begin{align}
\label{kQmin2}
&{}    (n^2-r)(n+r - n^2) +1 \leq  k  \leq (n^2-r) n , 
 \qquad  Q   \geq  \frac{n}{ n^2 - r } . 
\end{align}
\end{propnn}

Let us remark that the lower bound for $k$ given by Proposition~\ref{Qmin} is better than 
the one given by Proposition~\ref{rangeK} if $r$ is sufficiently small ($r \lesssim n/2$) 
or sufficiently large ($r \gtrsim n^2 -n/2$).

Recall that $\tau(\SR_i)$ is unitary if $Q \leq 2$. In the special case, when 
$Q=2$, equation (\ref{defH1}) acquires the form $\SR_i^2 = {\mathsf 1}$ 
and so $\tau(\SR_i)$ is an involutive solution to the Yang-Baxter equation. 
A classification of such solutions in terms of Thoma parameters was given in~\cite{LPW}. 
In particular, it was shown there that, for $Q=2$, a Temperley-Lieb type 
solution $P$ exists if and only if $\sqrt{n^2-4r}$ is an integer. 
Proposition~\ref{Qmin} imposes certain restrictions on the parameters 
of a unitary solution and an involutive solution in the general case.  

\begin{propnn}\label{Q2ns}
For $Q \leq 2$ and $n \geq 2r $, there exist no nontrivial solutions to (\ref{PPPm}),   
except in the case of $Q=2$, $n=2$, $r=1$. 
For $Q \leq 2$ and $2r \geq 2n^2 -n$,  there exist no nontrivial solutions to (\ref{PPPm}), 
except in the case of $Q=2$, $n=2$, $r=3$. 
\end{propnn}

\section{Examples of solutions not of the Temperley-Lieb type}

Examples of explicit Temperley-Lieb type solutions $P$ for all $n$, $r$ such that 
$n \geq r$ and $r \leq 4$ were given in \cite{By3} 
(cf., in particular, Theorem 2). Here we will consider an approach 
to constructing solutions to (\ref{PPPm}) in the general case and 
will provide some examples.

Let $G$ be the following  functional on   
$\mathrm {Mat}(r n,{\mathbb C}) \,{\times}\, {\mathbb N}$
\begin{align}
\label{Gak} 
{}&  G[A;k] = \bigl(k - r n + \tr A \bigr)^2 - k \, \big(k - r n  + \tr A^2 \bigr) .
\end{align} 

\begin{lemma}\label{Gcon}
Let $P_\CT $ be given by (\ref{PiTau}),
where $V_1, \ldots, V_r$ satisfy condition~(\ref{vv}), and let 
$A_\CT$ be given by~(\ref{WVAV}). If $P_\CT$ is a solution 
to (\ref{PPPm}), then there exists $k \in \mathbb N$ 
such that  $G[A_\CT;k] =0$. 
\end{lemma} 

To give an example of an application of Proposition \ref{AEsol} and Lemma \ref{Gcon}, 
let us consider the case $n=r=2$ and take $V_1$, $V_2$ of the following form
\begin{align}
\label{Vn2r2} 
{}&  V_1 = \frac{1}{\sqrt{a^2+1}} 
 \biggl( \begin{matrix}
  a e^{i \alpha} & 0 \\   0  & 1
 \end{matrix} \biggr) , \quad 
 V_2 = \frac{1}{\sqrt{b^2+1}} 
 \biggl( \begin{matrix}
  0 & 1   \\  b e^{i \beta}  & 0
 \end{matrix} \biggr)   , \quad a, b, \alpha, \beta \in {\mathbb R}.
\end{align} 
A direct computation shows that for such $V_1$, $V_2$, the functional 
$G[A_\CT;k]$ acquires a simple form only for $k=2$, namely, 
$G[A_\CT;2] = 8a^2 (1 + b^4 - 2b^2 \cos 2\beta)(1+a^2)^{-2}(1+b^2)^{-2}$. 
The necessary condition $G[A_\CT;2] =0$ is satisfied if $a=0$ or 
if $\beta =0$, $b= \pm 1$. It is easy to check that if $a=0$, then 
the only eigenvalues of the matrix $A_\CT$ are $1$ and  $(b+b^{-1})^{-2}$.  
If $\beta =0$ and $b= \pm 1$, then the only eigenvalues of 
$A_\CT$ are $1$ and  $\frac{1}{4}(a^{2}-1)^2 (a^{2}+1)^{-2}$.  
Therefore, by Proposition~\ref{AEsol}, in both these cases, $P_\CT$ 
is a solution of rank $r=2$ to (\ref{PPPm}) not of the Temperley-Lieb type. 
The corresponding values of $Q$ are $Q=b+b^{-1}$ in the first case and $Q=q+q^{-1}$, 
where $q^2 =\bigl((a+1)/(a-1)\bigr)^{\pm 2}$, in the second case.  
The corresponding matrices $R =q I_n^{\otimes 2} - QP$ have the following form 
(in the second case, we have taken $b=-1$, $a= (1-q)/(1+q)$ and set $Q_\pm = q \pm q^{-1}$)
\begin{align}
\label{Rn2r2} 
R = 
\begin{pmatrix}
  b & 0 & 0 & 0 \\
  0 & b-b^{-1} & -e^{- i \beta} & 0 \\
  0 & -e^{i \beta} & 0 & 0 \\
  0 & 0 & 0 & -b^{-1}
 \end{pmatrix}  , \quad 
 R =  \frac{1}{2}
 \begin{pmatrix}
  Q_- +2 & 0 & 0 & e^{i \alpha} Q_- \\
  0 & Q_- & Q_+  & 0 \\
  0 & Q_+ & Q_- & 0 \\
  e^{-i \alpha} Q_- & 0 & 0 & Q_- -2 
\end{pmatrix} . 
\end{align} 
The first $R$-matrix in (\ref{Rn2r2}) is related to the quantum supergroup $GL_q(1|1)$, see~\cite{ChKu},
whereas the second falls in the class of $R$-matrices satisfying the free fermion condition, 
$R_{11} R_{44} + R_{23} R_{32}  = R_{14} R_{41} + R_{22} R_{33}$. 

Now let us construct two examples for $n=3$ which appear to be novel. 
First, we take $V_1$, $V_2$, $V_3$ of the following form 
(where $a, b, c , \alpha, \beta \in {\mathbb R}$ and $d = (a^2 +b^2+c^2)^{\frac{1}{2} }$)
\begin{align}
\label{Rn3r3} 
 V_1 =  \frac{1}{d} 
\begin{pmatrix}
  0 & a & 0   \\
  b & 0 & 0 \\
  0 & 0 & c e^{i \alpha} 
 \end{pmatrix}  , \quad  
 V_2= \frac{1}{d} 
\begin{pmatrix}
  0 & 0 & b   \\
  0 & c e^{-i (\alpha+\beta)}  & 0 \\
  a & 0 & 0 
 \end{pmatrix}  , \quad  
 V_3 = \frac{1}{d}  
 \begin{pmatrix}
  c e^{i \beta} & 0 & 0   \\
  0 & 0 & a \\
  0 & b & 0 
 \end{pmatrix}  . 
\end{align} 
Then the functional $G[A_\CT;k]$ is a ratio of two polynomials in $a, b, c$ of degree 8  
and its numerator factorizes in a simple way only for $k=8$, namely, 
$G[A_\CT;8] = d^{-8} (a+b+c)^2 S$, where $S$ is a polynomial of degree~6. 
The necessary condition $G[A_\CT;8] =0$ is satisfied if  
\begin{align}
\label{abc0} 
 a+b+c =0 .
\end{align} 
It is easy to check that if the condition (\ref{abc0}) holds, then  
$A_\CT = \frac{1}{4} I_9 +  \frac{1}{4} x^t \,{\otimes}\, x$, 
where $x=(0,0,1,0,1,0,1,0,0)$. Since the only eigenvalues of $A_\CT$ 
are $1$ and $\frac{1}{4}$, Proposition~\ref{AEsol} implies that $P_\CT$ 
is a solution of rank $r=3$ to equations (\ref{PPPm}) for $Q=2$. 

In order to construct the second example, we take 
\begin{equation}
\begin{aligned}
\label{Rn3r4} 
{}&
 V_1 =  \frac{1}{\sqrt{|z_1|^2+|z_2|^2}}
\begin{pmatrix}
  0 & z_1 & 0   \\
  z_2 & 0 & 0 \\
  0 & 0 & 0 
 \end{pmatrix}  , \quad  
 V_2 =  \frac{1}{\sqrt{|z_1|^2+|z_2|^2}}
 \begin{pmatrix}
  0 & 0 & 0   \\
  0 & 0 & z_2 \\
  0 & z_1 & 0 
 \end{pmatrix} , \\ 
{}& 
 V_3=  \frac{1}{\sqrt{|z_3|^2+|z_4|^2}}
\begin{pmatrix}
  z_3 & 0 & 0   \\
  0 & 0  & 0 \\
 0 & 0 & z_4 
 \end{pmatrix}  , \quad  
 V_4 = \frac{1}{\sqrt{2}}  
 \begin{pmatrix}
  0 & 0 & 1   \\
  0 & 0 & 0 \\
  1 & 0 & 0 
 \end{pmatrix} ,
\end{aligned} 
\end{equation}
where $z_1, z_2 \in {\mathbb C}^*$, $z_3, z_4 \in {\mathbb C}$. 
Then the functional $G[A_\CT;k]$ is given by a ratio of two polynomials 
in $|z_1|$, $|z_2|$, $|z_3|$, $|z_4|$  
and its numerator factorizes for $k=2,6,8$. However, for $k=2,6$, 
the appearing factors have no real roots. In the remaining case, the numerator factorizes 
into four simple factors and 
the necessary condition $G[A_\CT;8] =0$ is satisfied if  
\begin{align}
\label{a1a2b1b2} 
  |z_3| (|z_1| - |z_2|) = \pm |z_4| (|z_1| + |z_2|)  \ \ \text{ or } \ \  
  |z_4| (|z_1| - |z_2|) = \pm |z_3| (|z_1| + |z_2|) .
\end{align} 
If the condition (\ref{a1a2b1b2}) holds, then 
the only eigenvalues of the matrix $A_\CT$ are $1$ and  $(q+q^{-1})^{-2}$, 
where $q=|z_1|/|z_2|$.  Therefore, 
$P_\CT$ is a solution of rank $r=4$ to equations (\ref{PPPm}) for  
$Q= |z_1|/|z_2| + |z_2|/|z_1|$. 

\appendix 
\renewcommand{\theequation}{A\arabic{equation}}
\setcounter{equation}{0}
\section{Appendix. Proofs.}

\begin{proof} [\bf Proof of Lemma~\ref{trPPrn}]   
a)  Since $(P_{1} - P_{2})^2$ is a positive semi-definite matrix, we have 
$\tr_3(P_{1} - P_{2})^2 = 2 r n -2 \tr_3 (P_1 P_2) \geq 0$.
The equality is possible only if $P_{1} = P_{2}$.  
But then we have $P_{2} = P_{3}$ as well. 
Hence $P \,{\otimes}\, I_n^{\otimes 2} = P_1= P_3 =  I_n^{\otimes 2} \,{\otimes}\, P$.
Taking the partial trace $\tr_2$ over the first factor in the Kronecker product, we get 
$r I_n^{\otimes 2} = n^2 P$, which implies that either $r=0$ and so $P=0$ or
$r=n^2$ and $P= I_n^{\otimes 2}$. \\[0mm]  
b) $H = i (P_1 P_2 - P_2 P_1)$ is a Hermitian matrix and so 
$H^2$ is positive semi-definite. Therefore, 
$\tr_3 H^2 = 2 \tr_3 (P_1 P_2) - 2 \tr_3 (P_1 P_2)^2 \geq 0$ 
and the equality is possible only if $H=0$.
\end{proof}

\begin{proof}[\bf Proof of Lemma~\ref{Fpos}] 
Given an orthogonal projection $P \in \mathrm {Mat}(n^2,{\mathbb C})$
of rank $r$, consider a one-parameter family of Hermitian matrices 
$H_\alpha \in \mathrm {Mat}(n^3,{\mathbb C})$, $\alpha \in {\mathbb R}$, where  
$H_\alpha  = P_{1} P_{2} P_{1} - P_{2} P_{1} P_{2} 
   - \alpha P_{1} +  \alpha P_{2}$. It is straightforward to check that
\begin{equation}\label{trHal}
\begin{aligned} 
{}&   \frac{1}{2} \tr_3 H_\alpha^2   =   a \, \alpha^2 + b \, \alpha +c , \quad  \text{where} \quad 
 a =  r n - \tr_3 (P_1 P_2) ,   \\  
{}& 
 b =  2 \tr_3 (P_1 P_2)^2 - 2 \tr_3 (P_1 P_2)   , \quad
 c =    \tr_3 (P_1 P_2)^2 - \tr_3 (P_1 P_2)^3  .   
\end{aligned}  
\end{equation}
If $a=0$, then by Lemma~\ref{trPPrn} we have $P=0$ or $P=I_n \,{\otimes}\, I_n$. 
In both these cases, we have $F[P]=0$ trivially. 
If $a \neq 0$, then we have 
$(a \, \alpha^2 + b \, \alpha +c) \geq 0$ for all $\alpha \in {\mathbb R}$  
since $H_\alpha^2$ is a positive semi-definite matrix. 
Consequently, the discriminant $\Delta = b^2 -4ac$ of this quadratic polynomial 
is non-positive.  It remains to notice that $F[P]=- \frac{1}{4} \Delta$. 
\end{proof}

\begin{proof}[\bf Proof of Proposition~\ref{tracecond}]
a)  If $P$ satisfies relations (\ref{PPPm}), then $\tr_3 (P_1 P_2)^m$ is given 
by~(\ref{trPPmk}) (see Proposition~\ref{TRsol}). 
Substituting the values given by this formula for $m=1,2,3$ in (\ref{trcond0}), 
we obtain 
$F[P]= k^2(1 -Q^{-2})(Q^{-4} -Q^{-6})-k^2 (Q^{-2} -Q^{-4})^2 =0$. \\[1mm]
b) Let $H_\alpha$ and $a,b,c, \Delta$ be the same as in the proof of Lemma~\ref{Fpos}. 
If  $P_1 P_2 \neq P_2 P_1$, then $P$ is not a multiple of $I_n \,{\otimes}\, I_n$ and 
hence, by Lemma~\ref{trPPrn}, $ a >0 $ and $b <0$. Therefore, we can set   
$\alpha_0 = - b/(2a)$ and we have $\alpha_0 >0$. 
Substituting $\alpha = \alpha_0$ in (\ref{trHal}), we obtain
\begin{align}
\nonumber
\frac{1}{2} \tr_3 H_{\alpha_0}^2   =   a \, \alpha_0^2 + b \, \alpha_0 +c 
= c - \frac{b^2}{4a} = -\frac{1}{4a} \Delta = \frac{1}{a} F[P]. 
\end{align}  
Since $H_{\alpha_0}^2$ is a positive semi-definite matrix, the condition $F[P]=0$
implies that $H_{\alpha_0}  = P_{1} P_{2} P_{1} - P_{2} P_{1} P_{2} 
   - \alpha_0 P_{1} +  \alpha_0 P_{2} =0$. That is,  
$P$ satisfies relations (\ref{PPPm}), where $Q^2= 1/{\alpha_0}$.
\end{proof}

\begin{proof}[\bf Proof of Lemma~\ref{dualP}]
An elementary computation.
\end{proof}

\begin{proof}[\bf Proof of Lemma~\ref{spekK}]
 We have 
$K_P^3 =   P_{1} - P_{2} - P_{1} P_{2} P_{1} + P_{2} P_{1} P_{2}$.
Therefore, the last relation in (\ref{PPPm})  is equivalent to the relation
$K_P^3 = (1 - Q^{-2}) K_P$. Since $K_P$ is Hermitian, we infer that
$P$ is a solution to (\ref{PPPm}) if and only if
the minimal polynomial of $K_P$ is $p(t)= t^3 + (Q^{-2} -1) t$. The roots  
of this polynomial are $0$ and $\lambda_\pm =  \pm \sqrt{1 - Q^{-2}}$. 
Since $K_P$ is Hermitian, $\lambda_\pm$ 
must be real, which implies that relations (\ref{PPPm}) do not admit 
nontrivial solutions if $Q \leq 1$.  
\end{proof}

\begin{proof}[\bf Proof of Proposition~\ref{rangeK}]
Since $\tr_3 K_P =0$, the eigenvalues 
$\lambda_+$  and $\lambda_-$ must have the same multiplicity~$k$. 
Hence $\mathrm{rank}\, K_P = 2k$.  
Therefore, we have $2k = \mathrm{rank}\, (P_1 - P_2) \leq 
  \mathrm{rank}\, (P \,{\otimes}\, I_n)  + \mathrm{rank}\, (I_n \,{\otimes}\, P) = 2 r n$. 
Using Lemma~\ref{dualP}, we also have 
$2k = \mathrm{rank}\, K_{\tilde{P}} \leq 
  \mathrm{rank}\, (\tilde{P} \,{\otimes}\, I_n)  + \mathrm{rank}\, (I_n \,{\otimes}\, \tilde{P}) = 
  2 \tilde{r}n = 2(n^2-r)n$.

In order to obtain a lower bound for $k$, we will use the following 
statement (see Theorem~1.3 in~\cite{ChTi}):  
if $X \in  \mathrm {Mat}(m,{\mathbb C})$ is an idempotent matrix, then 
$\mathrm{rank}\, (X \,{\otimes}\, I_m - I_m  \,{\otimes}\, X)  =2m \,\mathrm{rank}\, X -
 2 (\mathrm{rank}\, X)^2$. Therefore, for $(P_1 - P_3) \in \mathrm {Mat}(n^4,{\mathbb C})$, 
we have 
$\mathrm{rank}\, (P_1 - P_3) = \mathrm{rank}\, (P \,{\otimes}\,  I_{n^2}
- I_{n^2} \,{\otimes}\, P)= 2 n^2 \mathrm{rank}\, P - 2 (\mathrm{rank}\, P)^2
= 2 r(n^2 -r)$. On the other hand, 
$\mathrm{rank}\, (P_1 - P_3) = \mathrm{rank}\, (P_1 -  P_2 + P_2- P_3) \leq 
\mathrm{rank}\, (P_1 -  P_2) + \mathrm{rank}\, (P_2 -  P_3) = 
\mathrm{rank}\, (K_P \,{\otimes}\,  I_n) + \mathrm{rank}\, (I_n \,{\otimes}\,  K_P)
= 4 k n$. Thus, $4kn \geq 2 r(n^2 -r)$. 
\end{proof}

\begin{proof}[\bf Proof of Proposition~\ref{TRsol}] 
Multiplying the last relation in (\ref{PPPm})  with $P_1$ or $(P_1 P_2)^m P_1$, $m \geq 1$, 
 from the left and taking the trace $\tr_3$, we obtain, respectively, the following equalities 
\begin{align}
\label{trace0a}
{}&   Q^2 \bigl(  \tr_3 (P_1 P_2) - \tr_3 (P_1 P_2)^2  \bigr)  
   = r n- \tr_3 (P_1 P_2) , \\ 
\label{trace0m}
{}&  Q^2 \bigl(  \tr_3 (P_1 P_2)^{m+1} - \tr_3 (P_1 P_2)^{m+2}  \bigr)    
 =  \tr_3 (P_1 P_2)^{m} - \tr_3 (P_1 P_2)^{m+1}  .
\end{align} 
Denote $a_0 = r n$ and $a_m = \tr_3 (P_{1} P_{2})^m$ for $m \geq 1$. Then 
equations (\ref{trace0a})--(\ref{trace0m}) acquire the form  
$ a_{m+1} - a_m = Q^{-2} (a_m - a_{m-1})$, $m = 1,2, \ldots$. 
For $m \geq 0$, set $b_m = a_{m+1} - a_{m}$. Then we have 
$b_m = Q^{-2m} \, b_0$, where $Q>1$ by Lemma~\ref{spekK}. Hence 
\begin{align}\label{amm}
  \tr_3 (P_{1} P_{2})^m =
  a_m = a_0 + \sum_{i=0}^{m-1} b_i 
   = r n + \frac{1- Q^{-2m}}{1- Q^{-2}} ( \tr_3 (P_{1} P_{2}) - r n) .
\end{align}    
Next we note that  $\tr_3 K_P^2 = \tr_3 (P_1 + P_2 - 2 P_1 P_2) = 2 ( r n - \tr_3 (P_1 P_2) ) $.
On the other hand, Lemma~\ref{spekK} implies that the only non-zero eigenvalue of $K_P^2$ is  
$(1- Q^{-2})$. Whence $\tr_3 K_P^2 = 2 k (1- Q^{-2})$.
Thus, we have $\tr_3 (P_1 P_2) - r n = k (Q^{-2} -1)$. Substituting this equality in (\ref{amm}), 
we obtain equation~(\ref{trPPmk}). \\[1mm]
b) Substituting in (\ref{trcond0}) the values of $\tr_3 (P_{1} P_{2})^m$ 
given by equation (\ref{trPPmk}) for $m=1,2,3$, we check that  
$F[P]=0$. Noticing that, for $Q > 1$, the r.h.s.~of (\ref{trPPmk}) nontrivially depends 
 on $m$, we infer that $P_1 P_2 \neq P_2 P_1$. 
So, invoking the part b) of Proposition~\ref{tracecond}, we conclude that 
$P$ is a solution to~ (\ref{PPPm}). 
\end{proof}  

\begin{proof}[\bf Proof of Proposition~\ref{nonTL}] 
a) Set $X = Q^2   P_{1} P_{2} P_{1} - P_{1}$. 
It is straightforward to check that $X^2  = (Q^2-1) X$.  
Taking into account that  $Q > 1$ by Lemma \ref{spekK}, we infer that  
 $  (Q^2-1)^{-1} X$ 
is an orthogonal projection. Hence, using equation (\ref{trPPmk}) for $m=1$,
we obtain 
\begin{align}
\nonumber 
{}& 
 \mathrm {rank }\, X = \tr_3 \bigl( (Q^2-1)^{-1} X \bigr) = 
 (Q^2-1)^{-1} \bigl( Q^2 \tr_3 (P_1 P_2) - r n \bigr) 
 \stackrel{(\ref{trPPmk}) }{=}  r n -k . 
\end{align}
The second equality in (\ref{rankPPP}) can be derived in the same way.\\[1mm]
b) It is clear that $P$ is of the Temperley-Lieb type if and only if $\mathrm {rank }\, X = 0$.
\end{proof}

\begin{proof}[\bf Proof of Lemma~\ref{trAP}] 
Using the property 
$E^{(n)}_{ab}   E^{(n)}_{cd}  = \delta_{bc} E^{(n)}_{ad}$ of matrix units, 
it is straightforward to check that 
\begin{align}
\label{PYL}
{}&  (P_\CT)_{1} (P_\CT)_{2} (P_\CT)_{1} = 
 \sum_{s_1, s_2=1}^r Y_{s_1 s_2}  \otimes   L_{s_1 s_2 } ,
\end{align}
where  $Y_{s_1 s_2} \in  \mathrm {Mat}(n^2,{\mathbb C})$ and 
$L_{s_1 s_2 } \in  \mathrm {Mat}(n,{\mathbb C})$ are given by  
\begin{align}
\label{YLdef}
{}& Y_{s_1 s_2} =  \sum_{a,b,c,d=1}^n  (V_{s_1})_{ab} \, 
 (\bar{V}_{s_2})_{cd} \  E^{(n)}_{ac} \otimes E^{(n)}_{bd}  , \qquad
  L_{s_1 s_2} =   \sum_{i=1}^r V^t_{i} V_{s_1}^* V_{s_2} \bar{V}_{i} . 
\end{align}
Using condition (\ref{vv}), we obtain 
\begin{align}
\label{YYa}
{}& Y_{s_1 s_2} Y_{s_3 s_4} =  \delta_{s_2 s_3} Y_{s_1 s_4} , \qquad 
 \tr  Y_{s_1 s_2} = \delta_{s_1 s_2} .
\end{align}
Therefore, we have (repeated indices imply the summation)
\begin{align} 
\nonumber 
 {}& \tr_3 \bigl( (P_\CT)_{1} (P_\CT)_{2} \bigr)^m    =
 \tr_3 \bigl( (P_\CT)_{1} (P_\CT)_{2}  (P_\CT)_{1} \bigr)^m   \\[0.5mm] 
\nonumber 
{}&  = \tr \bigl( Y_{s_1 s_1'} Y_{s_2 s_2'} \ldots Y_{s_m s_m'} \otimes
  L_{s_1 s_1'} L_{s_2 s_2'} \ldots L_{s_m s_m'} \bigr) \\ 
\nonumber 
{}&  = \tr \bigl( Y_{s_1 s_1'} Y_{s_2 s_2'} \ldots Y_{s_m s_m'} \bigr) \, 
  \tr \bigl( L_{s_1 s_1'} L_{s_2 s_2'} \ldots L_{s_m s_m'} \bigr)   
  \stackrel{(\ref{YYa}) }{=}
 \tr \bigl( L_{s_1 s_2} L_{s_2 s_3} \ldots L_{s_m s_1} \bigr) . 
\end{align}
On the other hand, we also have (repeated indices again imply the summation)
\begin{align} 
\nonumber 
{}&  \tr (A_\CT)^m   =  
  \tr \bigl( E^{(r)}_{i_1 j_1} E^{(r)}_{i_2 j_2} \ldots E^{(r)}_{i_m j_m} 
  \otimes 
   (   V_{s_1} \bar{V}_{i_1} V^t_{j_1} V_{s_1}^* )  
   (   V_{s_2} \bar{V}_{i_2} V^t_{j_2} V_{s_2}^* )  \ldots 
    (   V_{s_m} \bar{V}_{i_m} V^t_{j_m} V_{s_m}^* )  \bigr)  \\[0.5mm]
\nonumber 
{}&  = \tr \bigl( E^{(r)}_{i_1 j_1} E^{(r)}_{i_2 j_2} \ldots E^{(r)}_{i_m j_m} \bigr)
 \tr \bigl( 
   (   V_{s_1} \bar{V}_{i_1} V^t_{j_1} V_{s_1}^* )  
   (   V_{s_2} \bar{V}_{i_2} V^t_{j_2} V_{s_2}^* )  \ldots 
   (   V_{s_m} \bar{V}_{i_m} V^t_{j_m} V_{s_m}^* )  \bigr) \\
\nonumber 
{}& =  \tr \bigl( 
   (   V_{s_1} \bar{V}_{i_1} V^t_{i_2} V_{s_1}^* )  
   (   V_{s_2} \bar{V}_{i_2} V^t_{i_3} V_{s_2}^* )  \ldots 
   (   V_{s_m} \bar{V}_{i_m} V^t_{i_1} V_{s_m}^* )  \bigr) 
 = \tr \bigl( L_{s_1 s_2} L_{s_2 s_3} \ldots L_{s_m s_1} \bigr) . 
\end{align} 
Thus, equality (\ref{trPPmkA}) is proved. 
\end{proof}

\begin{proof}[\bf Proof of Lemma~\ref{QNrk1}] 
For a trivial solution, we have either $r=k=0$ or $r=n^2$, $k=0$. In both 
these cases, relation (\ref{Qrn0}) holds trivially.

Let $P_\CT$ be a nontrivial solution to (\ref{PPPm}) 
and $W_\CT$, $A_\CT$ be given by~(\ref{WVAV}).  
Let $\|X\|_1 = \tr \bigl(X X^*)^{\frac 12}$ denote the trace norm and  
$\|X\|_2 = \bigl(\tr (X X^*) \bigr)^{\frac 12}$ denote the Frobenius norm.
Using the relation 
$\|X Y\|_1 \leq \|X\|_2 \|Y\|_2$ (see Proposition~9.3.6 in \cite{Ber}), we obtain 
\begin{align}
\nonumber
{}&  \| W_\CT \|_1 \leq \sum_{s,m=1}^r \| V_m \bar{V}_s \|_1 
 \leq \sum_{s,m=1}^r \| V_m \|_2 \| V_s \|_2  
 \stackrel{(\ref{vv}) }{=} r^2 .
\end{align}
On the other hand, equation (\ref{trAmk}) implies that 
$(A_\CT)^{\frac 12}$ has eigenvalue 1 with multiplicity $(r n -k)$ and 
eigenvalue $Q^{-1}$ with multiplicity~$k$. Hence, 
$r^2  \geq \| W_\CT \|_1  = \tr (A_\CT)^{\frac 12} = r n - k + k/Q$.
\end{proof}

\begin{proof}[\bf Proof of Proposition~\ref{TLrank1}] 
Let $r=1$. Then $k \neq 0$ since $P$ is a nontrivial solution 
and we have $k \leq n$ by Proposition~\ref{rangeK}. 
On the other hand, equation (\ref{Qrn0}) implies that $n-k \leq 1 - k/Q <1$. 
Thus, we conclude that $k=n$ which by Proposition~\ref{nonTL}b 
implies that $P$ is of the Temperley-Lieb type.
\end{proof}

\begin{proof}[\bf Proof of Proposition~\ref{Qmin}]  
a) Inequality (\ref{Qrn0}) implies that $r n  - r^2 \leq k - k/Q < k$. Hence 
$k \geq r n -r^2 +1$. Since $r < n \leq n^2/2$, we have $k \leq r n$ by 
Proposition~\ref{rangeK}. Moreover,  by Proposition \ref{nonTL}b, we have 
$k \neq r n$ since $P$  is not of the Temperley-Lieb type. 
Thus, $k  \leq r n - 1$. 
Inequality (\ref{Qrn0}) also implies that $(1-Q^{-1}) \geq (r n - r^2)/k$.  Hence, 
 $(1-Q^{-1}) \geq (r n - r^2)/(r n -1)$, 
which is equivalent to the second inequality in~(\ref{kQmin1}). \\[1mm]
b) The dual solution $\tilde{P}$ has the same value of~$k$ and its 
rank is $\tilde{r} = n^2 - r < n$ so that we have 
$k \geq \tilde{r} n - \tilde{r}^2 +1$ by the part a) of this proposition. 
Since $r > n^2 - n \geq n^2/2$, we have $k \leq \tilde{r} n$ by 
Proposition~\ref{rangeK}. Note that, unlike in the part a), we cannot claim  
that $k \neq \tilde{r} n$ because we cannot exclude apriori that the 
dual solution is of the Temperley-Lieb type. 
 The dual solution $\tilde{P}$ has the same value of~$Q$ and 
inequality (\ref{Qrn0}) implies that $(1-Q^{-1}) \geq (\tilde{r} n - \tilde{r}^2)/k$.
Hence $(1-Q^{-1}) \geq (\tilde{r} n - \tilde{r}^2)/(\tilde{r} n )$, 
which is equivalent to the second inequality in~(\ref{kQmin2}).
\end{proof}

\begin{proof}[\bf Proof of Proposition~\ref{Q2ns}] 
Let $P \in \mathrm {Mat}(n^2,{\mathbb C})$ be a nontrivial solution to (\ref{PPPm})  
for $n \geq 2r$. If $P$ is not of the Temperley-Lieb type, then 
the second inequality in (\ref{kQmin1}) yields $Q \geq (2r^2-1)/(r^2-1) >2$. 
If $P$ is of the Temperley-Lieb type and $n>2r$,  then $Q \geq n/r >2$ (by Theorem 3 in \cite{By1} 
or as a corollary to Lemma~\ref{QNrk1}). There remains the case where 
$P$ is of the Temperley-Lieb type and $Q \leq 2$, $n=2r$. However, for $Q=2$, a solution to 
(\ref{TTTm}) exists if and only if $\sqrt{n^2-4r}$ is an integer (see Proposition 6.3 in \cite{LPW}). 
But $\sqrt{(2r)^2-4r} = 2 \sqrt{r(r-1)}$ is not an integer unless $r=1$.  
Finally, for $Q < 2$, a nontrivial solution to (\ref{TTTm}) exists only if  $n^2 =2r$ or 
$n^2=3r$ (see Theorem~1 in~\cite{By3}) which is not consistent with 
the condition $n \geq 2r$ if $n \geq 2$.\\ 
Let $P \in \mathrm {Mat}(n^2,{\mathbb C})$ be a solution to (\ref{PPPm})  
for $Q \leq 2$ such that $2r \geq 2n^2 -n$. 
Then for the dual solution $\tilde{P}$, we have the same value of $Q$ and  
$n \geq 2 \tilde{r}$.  By the part a) of this Proposition, this is possible 
only if $Q=2$, $n=2$, and $\tilde{r}=1$, that is $r=3$.
\end{proof}

\begin{proof}[\bf Proof of Lemma~\ref{Gcon}]  
Equations (\ref{trAmk}) for $m=1,2$ read  
$(k - r n + \tr A_\CT ) = Q^{-2} k$ and 
$(k - r n + \tr A_\CT^2 ) = Q^{-4} k$. Excluding the parameter $Q$, 
we obtain an equality equivalent to the condition $G[A_\CT;k] =0$. 
\end{proof}

\vspace*{1mm}
\noindent 
\small{
{\bf Acknowledgements.} 
This work was supported in part by the program NCCR SwissMAP of
the Swiss National Science Foundation.   
The author is thankful to the referee, who reviewed an earlier version 
of this manuscript submitted for publication, for useful comments. 
}

 \vspace*{1.5mm}
{\sc \small
\noindent
Section of Mathematics, University of Geneva,  
C.P. 64, 1211 Gen\`eve 4, Switzerland  \\[1mm] 
Steklov Mathematical Institute, 
Fontanka 27, 191023, St. Petersburg, Russia  
}

\end{document}